\documentclass[journal]{IEEEtran}
 
\usepackage[T1]{fontenc}
\usepackage[utf8]{inputenc}
 
\usepackage{cite}      

\usepackage[dvips]{graphicx}  

\usepackage{amsmath,amsfonts,amssymb}

\DeclareMathOperator*{\minimize}{minimize}

\begin{document}

\title{Tikhonov Regularization of Sphere-Valued Signals}

\author{Laurent~Condat\\[3mm]
July 25, 2022
\thanks{L. Condat is with the Visual Computing Center, King Abdullah University of Science and Technology (KAUST), Thuwal 23955-6900, Kingdom of Saudi Arabia.
Contact: see https://lcondat.github.io/}
}

\markboth{Tikhonov Regularization of Sphere-Valued Signals}{}
\maketitle

\begin{abstract}
It is common to have to process signals, whose values are points on the 3-D sphere. We consider a Tikhonov-type regularization model to smoothen or interpolate sphere-valued signals defined on arbitrary graphs. We propose a convex relaxation of this nonconvex problem as a semidefinite program, which is easy to solve numerically and is efficient in practice. 
\end{abstract}

\IEEEpeerreviewmaketitle

\section{Introduction}

\IEEEPARstart{I}{n} a wide range of applications, one has to deal with signals or data records, like time series, whose values are points on the 3-D unit sphere, like the orientation of a camera, direction of a moving object, axis of rotation of a robot's articulation. 
We consider the general setting, where the (possibly multidimensional) signal is defined on a graph, with values located at the nodes; every value is a point on the sphere. Two values are adjacent if there is an edge between their nodes. A 2-D image is a particular case with edges between every pair of neighboring pixels horizontally and vertically, forming a square grid. Then, to regularize signals on graphs, it is natural to promote the property that adjacent values are close to each other, in some sense. For scalar values, Tikhonov regularization consists in penalizing the squared differences of adjacent values. In this work, we focus on an equivalent model for sphere-valued signals. The sphere is nonconvex, so that the corresponding Tikhonov regularization problem is nonconvex too, and challenging to solve. We propose a new convex relaxation and show by some numerical examples that it is tight enough to yield exactly the \emph{global} solution of the Tikhonov nonconvex problem in practical cases.

\section{Tikhonov Smoothing for Sphere-Valued Signals}\label{sec2}

\subsection{Sphere-Valued Signals on Graphs}

We denote by $\mathbb{S}=\{x=(x^1,x^2,x^3)\in\mathbb{R}^3\ :\ \|x\|_2=1\}$ the 3-D unit sphere, where $\|\cdot\|_2$ is the Euclidean norm.
We want to estimate a signal $x=(x_n)_{n\in V}$, with values $x_n \in \mathbb{S}$, 
defined on a connected undirected graph $(V,E)$,  
where 
$V$
is the set of nodes 
and $E$ is the set of edges, which are sets of two distinct nodes.
Typically, we are given a noisy signal $y=(y_n)_{n\in V}$ defined on the same graph and the sought signal $x$ is a smoothed, or denoised, version of $y$, which achieves a tradeoff between closeness to $y$ and smoothness, in some sense. Another typical setting is interpolation, or inpainting: $y$ is defined on a subset $U\subset V$ of nodes and we want to estimate its missing samples; that is, $x$ is the smoothest signal defined on $V$ such that $x_n=y_n$, for every $n\in U$.

\subsection{Classical Tikhonov Regularization}  

For vector-valued signals with values in $\mathbb{R}^d$, for any $d\geq 1$, Tikhonov-regularized smoothing consists in solving the following convex optimization problem. 
Given $y=(y_n)_{n\in V}$ and nonnegative weights $(w_n)_{n\in V}$ and $(\lambda_{n,n'})_{\{n,n'\}\in E}$, $x=(x_n)_{n\in V}$ is the solution to 
\begin{equation}
\minimize_{x_n \in \mathbb{R}^d\,:\, n\in V}
 \,\sum_{n\in V} \frac{w_n}{2} \|x_n-y_n\|_2^2 + \! \!\sum_{\{n,n'\}\in E}\! \frac{\lambda_{n,n'}}{2} \|x_n-x_{n'}\|_2^2.\label{eq1}
\end{equation}
For the interpolating task with $y$ defined only on $U\subset V$, we want to solve, instead:
\begin{equation}
\minimize_{x_n \in \mathbb{R}^d\,:\, n\in V}  \!\sum_{\{n,n'\}\in E}\! \frac{\lambda_{n,n'}}{2} \|x_n-x_{n'}\|^2\quad\mbox{s.t.}\quad x_n=y_n,\ \forall n\in U.\label{eq2}
\end{equation}
Formally, \eqref{eq2} can be viewed as a particular case of \eqref{eq1} with $w_n=\{+\infty$ if $n\in U$,  0 otherwise$\}$, so that we can focus on the form  \eqref{eq1}, with the weights $w_n$ allowed to be $+\infty$.

We can notice that \eqref{eq1} corresponds to the maximum-a-posteriori (MAP) estimate of an unknown signal $x^\sharp$ given $y$, which is $x^\sharp$ plus white Gaussian noise, assuming a Gaussian Markov Random Field  prior for $x^\sharp$, with nonzero dependencies between its Gaussian variables along the edges of $V$. That is, $y_n^k -x^{\sharp, k}_n\sim \mathcal{N}(1/w_n)$ and  $x^{\sharp,k}_n-x^{\sharp,k}_{n'} \sim \mathcal{N}(1/\lambda_{n,n'})$, where $\mathcal{N}(\sigma^2)$ denotes the normal distribution with zero mean and variance $\sigma^2$, independently on each coordinate $k=1,\ldots,d$.

\subsection{Proposed Model}

We want to formulate an equivalent problem to \eqref{eq1} for signals $x$ and $y$ with values in $\mathbb{S}$. The natural counterpart of  \eqref{eq1}  is to keep the same cost function to minimize, with values in $\mathbb{S}$ instead of $\mathbb{R}^d$:
\begin{equation}
\minimize_{x_n \in \mathbb{S}\,:\, n\in V}
 \,\sum_{n\in V} \frac{w_n}{2} \|x_n-y_n\|_2^2 + \! \!\sum_{\{n,n'\}\in E}\! \frac{\lambda_{n,n'}}{2} \|x_n-x_{n'}\|_2^2.\label{eq1n}
\end{equation}
We can develop the squared norm $\|x_1-x_2\|^2$, for any $(x_1,x_2)\in\mathbb{S}^2$, as  $\|x_1\|^2+\|x_2\|^2-2 x_1\cdot x_2$, where $x_1\cdot x_2=x_1^1x_2^1 + x_1^2 x_2^2+x_1^3 x_2^3$ denotes the Euclidean inner product of $x_1$ and $x_2$. Since $\|x_1\|=\|x_2\|=1$, we have $\|x_1-x_2\|^2 = 2-2 x_1\cdot x_2$. Hence, we can rewrite \eqref{eq1n} as
\begin{align}
&\minimize_{x_n \in \mathbb{S}\,:\, n\in V}
 \ \Psi_\mathrm{orig}(x)=
 \sum_{n\in V}w_n(1-  x_n\cdot y_n)\notag \\
 &\quad+ \! \!\sum_{\{n,n'\}\in E}\! \lambda_{n,n'} (1-x_n \cdot x_{n'}).\label{eq1n2}
\end{align}
We can note that $x_n \cdot x_{n'} \in [-1,1]$ is the cosine of the angle at the origin between $x$ and $x'$, so that $1-x_n \cdot x_{n'}$ is indeed a kind of distance between $x_n$ and $x_{n'}$. 
Like in the classical setting, there is a Bayesian interpretation of \eqref{eq1n2} as a MAP estimate  of an unknown signal $x^\sharp$ given $y$, which is $x^\sharp$ plus noise, this time following a von Mises--Fisher distribution, 
widely used in directional statistics \cite{jup09}. That is, $y_n\in\mathbb{S}$ is the outcome of a random variable with probability density function (p.d.f.) proportional to 
$e^{w_n y_n \cdot x^\sharp_n}$. Similarly, $x^\sharp$ is supposed to be a Markov Random Field with von Mises--Fisher dependency along every edge, with p.d.f.\ proportional to $e^{\lambda_{n,n'} x^\sharp_{n'} \cdot x^\sharp_n}$.

\subsection{A Basic Convex Relaxation}

Because $\mathbb{S}$ is nonconvex, the problem \eqref{eq1n} is nonconvex and difficult to solve. A straightforward way to relax it into a convex problem is to minimize the objective function, which is convex, with every variable $x_n$ in the convex Euclidean ball $\mathbb{B}:=\{x\in\mathbb{R}^3\ :\ \|x\|_2\leq 1\}$ instead of the sphere $\mathbb{S}$. After this problem has been  solved, every $x_n \in \mathbb{B}$ is rescaled as $x_n / \|x_n\|$ to project it back on $\mathbb{S}$. We call the combination of this convex problem and the rescaling postprocessing step the baseline method.

\section{A New Convex Relaxation}

To formulate a better convex relaxation than simply replacing $\mathbb{S}$ by $\mathbb{B}$ in \eqref{eq1n}, we start with the equivalent problem \eqref{eq1n2} and introduce an auxiliary variable $d_{n,n'}\in\mathbb{R}$ at every edge $\{n,n'\}\in E$: we can rewrite \eqref{eq1n2}  as 
\begin{align}
&\minimize_{\small\substack{x_n \in \mathbb{S}\;:\; n\in V\\ d_{n,n'}\in \mathbb{R}\;:\; \{n,n'\}\in E}}
 \ \sum_{n\in V}w_n(1-  x_n\cdot y_n) \notag \\
 &\qquad+\sum_{\{n,n'\}\in E}\! \lambda_{n,n'} (1- d_{n,n'})\notag\\ 
&\qquad\mbox{s.t.}\ \ d_{n,n'}=x_n \cdot x_{n'},\ \forall  \{n,n'\}\in E.\label{eq1n3}
\end{align}
In this last problem, the objective function to minimize is convex, and even linear, which is beneficial. Indeed, a general property in optimization is that minimizing a linear functional over a set yields a solution on the boundary of the set.  It seems intractable to express the convex hull of the $x_n$ and $d_{n,n'}$ such that $x_n \in \mathbb{S}$ and  $d_{n,n'}=x_n \cdot x_{n'}$, for every $n,n'$. So, we ``marginalize'' the relaxation and design instead a convex set $\Omega_{n,n'}$ 
at every edge $\{n,n'\}\in E$, so that its boundary consists of points satisfying $x_n \in \mathbb{S}$, $x_{n'}  \in \mathbb{S}$ and  $d_{n,n'}=x_n \cdot x_{n'}$.

In a previous work \cite{con22}, the author proposed a convex relaxation of \eqref{eq1n3} for Tikhonov regularization on the 2-D unit circle $\{x=(x^1,x^2)\in\mathbb{R}^2\ :\ \|x\|_2=1\}$, where the $\Omega_{n,n'}$ are complex elliptopes, which are complex Hermitian matrices with particular structure. Indeed, it is convenient to represent points on the circle as complex numbers with modulus one, since multiplication of two such numbers corresponds to a rotation, or geodesic path between two points on the circle. Extending these ideas from the 2-D circle to the 3-D sphere is not straightforward, since there is no easy representation of 3-D rotations using complex numbers. However, it is known that quaternions are well suited to perform geometric operations on the sphere~\cite{arn95}, and that a vector quaternion of the form $q=a\mathrm{i} +b\mathrm{j} + c\mathrm{k}$, for some $(a,b,c)\in\mathbb{R}^3$, can be represented as a $2\times 2$ complex matrix, which is a scaled Pauli matrix:
\begin{equation}
M_{a,b,c}:=\left(\begin{array}{cc}
-ci&-b-ai\\
b-ai&ci
\end{array}
\right)
\end{equation}
(where the complex $i=\sqrt{-1}$ is different from the quaternionic $\mathrm{i}$ mentioned above).

For every $(a,b,c)\in\mathbb{R}^3$, $M_{a,b,c}^\mathrm{H} =-M_{a,b,c}$, where $\cdot^\mathrm{H}$ denotes the Hermitian transpose, and the product $M_{a,b,c}^\mathrm{H}M_{a,b,c}$ is 
\begin{align}
&\left(\begin{array}{cc}
ci&b+ai\\
-b+ai&-ci
\end{array}
\right)
\left(\begin{array}{cc}
-ci&-b-ai\\
b-ai&ci
\end{array}
\right)\\
&=\left(\begin{array}{cc}
a^2+b^2+c^2&0\\
0&a^2+b^2+c^2
\end{array}
\right)=(a^2+b^2+c^2) \mathrm{Id},\notag
\end{align}
where $\mathrm{Id}$ denotes the $2\times 2$ identity. Hence, $M_{a,b,c}^\mathrm{H}M_{a,b,c}=\mathrm{Id}$ if and only if $a^2+b^2+c^2=1$. More generally, for any $x=(a,b,c)\in\mathbb{S}$ and $x'=(a',b',c')\in\mathbb{S}$, the product $M_{a,b,c}^\mathrm{H}M_{a',b',c'}$ encodes information about the geodesic path 
between $x$ and $x'$:  for every $x=(a,b,c)\in\mathbb{R}^3$ and $x'=(a',b',c')\in\mathbb{R}^3$, the product $M_{a,b,c}^\mathrm{H}M_{a',b',c'}$ is
\begin{align}
&\left(\begin{array}{cc}
ci&b+ai\\
-b+ai&-ci
\end{array}
\right)
\left(\begin{array}{cc}
-c'i&-b'-a'i\\
b'-a'i&c'i
\end{array}
\right)\\
&\small=\left(\!\!\begin{array}{cc}
aa'\!+\!bb'\!+\!cc'+(ab'\!-\!a'b)i&a'c - ac'+(bc'\!-\!b'c)i\\
ac'-a'c+(bc'\!-\!b'c)i&aa'\!+bb'\!+cc'+(a'b\!-\!ab')i
\end{array}\!\!
\right)\!.\notag
\end{align}
The trace of this matrix equals twice the inner product $x\cdot x'= aa'+bb'+cc'$, which we are interested in for our formulation \eqref{eq1n3}. 

Therefore, for every $x=(a,b,c)\in\mathbb{R}^3$ and $x'=(a',b',c')\in\mathbb{R}^3$, the Hermitian matrix
\begin{equation}
\small\left(\begin{array}{cc|cc|cc}
1&0&-ci&-b-ai&-c'i&-b'-a'i\\
0&1&b-ai&ci&b'-a'i&c'i\\
\hline
ci&b+ai&1&0&d-gi&-f-ei\\
-b+ai&-ci&0&1&f-ei&d+gi\\
\hline
c'i&b'+a'i&d+gi&f+ei&1&0\\
-b'+a'i&-c'i&-f+ei&d-gi&0&1
\end{array}
\right)\label{eqmat1}
\end{equation}
is of rank 2 and positive semidefinite, in which case it can be written as the product
\begin{align}
&\left(\begin{array}{cc}
1&0\\
0&1\\
\hline
ci&b+ai\\
-b+ai&-ci\\
\hline
c'i&b'+a'i\\
-b'+a'i&-c'i
\end{array}
\right)\\
&\times \left(\begin{array}{cc|cc|cc}
1&0&-ci&-b-ai&-c'i&-b'-a'i\\
0&1&b-ai&ci&b'-a'i&c'i
\end{array}
\right),\notag
\end{align}
if and only if $a^2+b^2+c^2=1$, ${a'}^2+{b'}^2+{c'}^2=1$, $d=aa'+bb'+cc'$, $e=b'c-bc'$, $f=ac'-a'c$, $g=a'b-ab'$.
Since the boundary of the convex set of positive semidefinite matrices of the form \eqref{eqmat1} is actually the subset of rank-2 matrices, such matrices are well suited to promote the properties that $x\in\mathbb{S}$, $x'\in\mathbb{S}$, $d=x \cdot x'$.

Hence, the proposed convex relaxation is:
\begin{align}
&\minimize_{\small\substack{x_n \in \mathbb{R}^3\;:\; n\in V\\ (d_{n,n'},e_{n,n'},f_{n,n'},g_{n,n'})\in \mathbb{R}^4\;:\; \{n,n'\}\in E}}
\hspace{-3mm}\Psi_\mathrm{conv}(x,d,e,f,g)\notag\\
 &\qquad=\ \sum_{n\in V}w_n(1-  x_n\cdot y_n) +\sum_{\{n,n'\}\in E}\! \lambda_{n,n'} (1- d_{n,n'})\notag\\ 
&\qquad\mbox{s.t.}\ \ P_{n,n'}\succcurlyeq 0,\ \forall  \{n,n'\}\in E,\label{eq1n4}
\end{align}
where $\succcurlyeq 0$ denotes semidefiniteness and, for every $\{n,n'\}\in E$, $P_{n,n'}$ is a Hermitian matrix of the form \eqref{eqmat1}, with $a$, $a'$, $b$, $b'$, $c$, $c'$, $d$, $e$, $f$, $g$ replaced by $x_n^1$,  $x_{n'}^1$, $x_n^2$,  $x_{n'}^2$, $x_n^3$,  $x_{n'}^3$, $d_{n,n'}$, $e_{n,n'}$, $f_{n,n'}$, $g_{n,n'}$, respectively. There is no guarantee that when solving the problem \eqref{eq1n4}, every matrix $P_{n,n'}$  will be of rank 2, but this is what we hope for.

Thus, we solve the convex problem \eqref{eq1n4}, and if the obtained solution $(x^\star,d^\star,e^\star,f^\star,g^\star)$ is such that $x_n \in \mathbb{S}$ fore every $n\in V $ and  $d_{n,n'}=x_n \cdot x_{n'}$, for every $\{n,n'\}\in E$, then $x^\star$ is the exact global solution to the original nonconvex problem \eqref{eq1n2}.

The problem \eqref{eq1n4} consists of minimizing a linear function under a composite constraint, like in \cite{con22}. A well suited algorithm is the Proximal Method of Multipliers \cite{roc762,con192}; its complexity is dominated by projecting, at every iteration and every edge, a $6\times 6$ matrix onto the cone of positive semidefinite matrices.

\section{Experiments}

In progress.

\nocite{ali10}

\bibliographystyle{IEEEtran}
\bibliography{IEEEabrv,../biblio.bib}

\end{document}